\documentclass[11pt,centertags,reqno]{amsart}
\usepackage{amsfonts}
\usepackage{amssymb}
\usepackage{graphicx}
\usepackage{amscd}

\theoremstyle{definition}

\theoremstyle{plain}
\numberwithin{equation}{section}
\input{tcilatex}

\begin{document}

\begin{center}
{\LARGE A shrinking projection approximant for the split equilibrium
problems and fixed point problems in Hilbert spaces\bigskip }

Abdul Ghaffar$^{1}$, Zafar Ullah$^{2},$ Muhammad Aqeel Ahmad Khan$^{3,\ast }$
and Faisal Mumtaz$^{3}$

$^{1}$Department of Mathematical Sciences, Balochistan University of
Information Technology, Engineering and Management Sciences, Quetta 87300,
Pakistan

$^{2}$Department of Mathematics, University of Education Lahore, DG Khan
Campus,

DG Khan 32200, Pakistan

$^{3}$Department of Mathematics, COMSATS Institute of Information Technology
Lahore,

Lahore, 54000, Pakistan

March 06, 2018 \let\thefootnote\relax\footnote{%
*Corresponding author\newline
E-mail addresses: (A. Ghaffar) abdul.ghaffar@buitms.edu.pk, (Z. Ullah)
zafarbhatti73@gmail.com, (M.A.A. Khan) itsakb@hotmail.com,
maqeelkhan@ciitlahore.edu.pk (F. Mumtaz) fmnikyana49@gmail.com}
\end{center}

\textbf{Abstract}: This work is devoted to establish the strong convergence
results of an iterative algorithm generated by the shrinking projection
method in Hilbert spaces. The proposed approximation sequence is used to
find a common element in the set of solutions of a finite family of split
equilibrium problems and the set of common fixed points of a finite family
of total asymptotically strict pseudo contractions in such setting. The
results presented in this paper improve and extend some recent corresponding
results in the literature.\\[1mm]

\section{Introduction and Preliminaries}

Throughout this paper, we write $x_{n}\rightarrow x~($resp. $%
x_{n}\rightharpoonup x)$ to indicate the strong convergence (resp. the weak
convergence) of a sequence $\{x_{n}\}_{n=1}^{\infty }$. Let $C$ be a
nonempty subset of a real Hilbert space $H$ and let $T:C\rightarrow C$ be a
mapping. The set of fixed points of the mapping $T$ is defined and denoted: $%
F(T)=\{x\in C:T(x)=x\}.$ A self-mapping $T$\ is said to be: (i) nonexpansive
if $\left\Vert Tx-Ty\right\Vert \leq \left\Vert x-y\right\Vert ;$ (ii)
asymptotically nonexpansive \cite{Goebel Kirk} if there exists a sequence $%
\{\lambda _{n}\}\subset \lbrack 0,\infty )$ with $\lim_{n\rightarrow \infty
}\lambda _{n}=0$ such that $\left\Vert T^{n}x-T^{n}y\right\Vert \leq
(1+\lambda _{n})\left\Vert x-y\right\Vert ,$ $n\geq 1;$ (iii) Lipschitzian
if $\left\Vert T^{n}x-T^{n}y\right\Vert \leq \Theta \left\Vert
x-y\right\Vert $ for some $\Theta >0;$ (iv) firmly nonexpansive if%
\begin{equation}
\left\Vert Tx-Ty\right\Vert ^{2}+\left\Vert (I-T)x-(I-T)y\right\Vert
^{2}\leq \left\Vert x-y\right\Vert ^{2};  \label{1.1}
\end{equation}%
(v) pseudo-contraction \cite{Browder 1967}, if%
\begin{equation}
\left\Vert Tx-Ty\right\Vert ^{2}\leq \left\Vert x-y\right\Vert
^{2}+\left\Vert (I-T)x-(I-T)y\right\Vert ^{2};  \label{1.2}
\end{equation}%
(vi) $k$-strict pseudo contraction \cite{Browder 1967}, if there exists $%
k\in \lbrack 0,1)$ such that 
\begin{equation}
\Vert Tx-Ty\Vert ^{2}\leq \Vert x-y\Vert ^{2}+k\Vert (I-T)x-(I-T)y\Vert ^{2};
\label{1.3}
\end{equation}%
(vii). $(k,\{\lambda _{n}\})$-asymptotically strict pseudo contraction \cite%
{Qihou}, if there exist a constant $k\in \lbrack 0,1)$ and a sequence $%
\{\lambda _{n}\}\subset \lbrack 1,\infty )$ with $\lim_{n\rightarrow \infty
}\lambda _{n}=1$ such that%
\begin{equation}
\Vert T^{n}x-T^{n}y\Vert ^{2}\leq \lambda _{n}\Vert x-y\Vert ^{2}+k\Vert
(I-T^{n})x-(I-T^{n})y)\Vert ^{2};  \label{1.4}
\end{equation}%
(viii). $(\{\lambda _{n}\},\{\mu _{n}\},\xi )$-total asymptotically
nonexpansive \cite{Alber FPTA} if there exist nonnegative real sequences $%
\left\{ \lambda _{n}\right\} _{n=1}^{\infty }$, $\left\{ \mu _{n}\right\}
_{n=1}^{\infty }$ with $\lim_{n\rightarrow \infty }\lambda
_{n}=0=\lim_{n\rightarrow \infty }\mu _{n}$ and a strictly increasing
continuous function $\xi :\mathbb{R^{+}}\rightarrow \mathbb{R^{+}}$ with $%
\xi (0)=0$ such that%
\begin{equation}
\Vert T^{n}x-T^{n}y\Vert \leq \Vert x-y\Vert +\lambda _{n}\xi (\Vert
x-y\Vert )+\mu _{n};  \label{1.5}
\end{equation}%
(ix). $(k,\{\lambda _{n}\},\{\mu _{n}\},\xi )$-total asymptotically strictly
pseudo contraction \cite{Yang-FPTA}, if there exist a constant $k\in \lbrack
0,1)$ and nonnegative real sequences $\left\{ \lambda _{n}\right\}
_{n=1}^{\infty }$, $\left\{ \mu _{n}\right\} _{n=1}^{\infty }$ with $%
\lim_{n\rightarrow \infty }\lambda _{n}=0=\lim_{n\rightarrow \infty }\mu
_{n} $ and a strictly increasing continuous function $\xi :\mathbb{R^{+}}%
\rightarrow \mathbb{R^{+}}$ with $\xi (0)=0$ such that%
\begin{equation}
\Vert T^{n}x-T^{n}y\Vert ^{2}\leq \Vert x-y\Vert ^{2}+k\Vert
(I-T^{n})x-(I-T^{n})y)\Vert ^{2}+\lambda _{n}\xi (\Vert x-y\Vert )+\mu _{n},
\label{1.6}
\end{equation}%
holds for all $x,y\in C.$\newline
\textbf{Remark 1.1.} It is worth mentioning that the class of nonexpansive
mappings have powerful applications to solve various problems arising in the
field of applied mathematics, such as variational inequality problem, convex
minimization, zeros of a monotone operator, initial value problems of
differential equations, game-theoretic model and image recovery. It is
therefore, natural to extend such powerful results of the class of
nonexpansive mappings to the more general class of mappings. As a
consequence, the notion of nonexpansive mapping has been generalized in
several ways. In 1967, Browder and Petryshyn \cite{Browder 1967} introduced
the concept of strict pseudo contraction as a generalization of nonexpansive
mappings. Later on, Alber et al.\cite{Alber FPTA} introduced the notion of
total asymptotically nonexpansive mappings which is more general in nature
and unifies various definitions of mappings associated with the class of
asymptotically nonexpansive mappings. In 2011, Yang et al. \cite{Yang-FPTA}
introduced the notion of total asymptotically strict pseudo contraction
which contains properly the class of total asymptotically nonexpansive
mappings and strict pseudo contractions. So, we study this general class of
mappings to contribute in metric fixed point theory.

Let $C$ be a nonempty subset of a real Hilbert space $H_{1}$, $Q$ be a
nonempty subset of a real Hilbert space $H_{2}$ and let $A:H_{1}\rightarrow
H_{2}$ be a bounded linear operator. Let $f:C\times C\rightarrow \mathbb{R}$
and $g:Q\times Q\rightarrow \mathbb{R}$ be two bifunctions. The split
equilibrium problem (SEP) is to find:%
\begin{equation}
x^{\ast }\in C\text{ \ \ \ such that }f\left( x^{\ast },x\right) \geq 0\text{
for all }x\in C,  \label{1.7}
\end{equation}%
and%
\begin{equation}
y^{\ast }=Ax^{\ast }\in Q\text{ \ \ \ such that }g\left( y^{\ast },y\right)
\geq 0\text{ for all }y\in Q.  \label{1.8}
\end{equation}%
It is remarked that inequality (1.7) represents the classical equilibrium
problem \cite{Combettes} and its solution set is denoted $EP(f).$ Moreover,
inequalities (1.7) and (1.8) constitute a pair of equilibrium problems which
aim to find a solution $x^{\ast }$\ of an equilibrium problem (1.7) such
that its image $y^{\ast }=Ax^{\ast }$ under a given bounded linear operator $%
A$ also solves another equilibrium problem (1.8). The set of solutions of
SEP (1.7) and (1.8) is denoted $\Omega =\{z\in EP(f):Az\in EP(g)\}.$

Equilibrium problem theory provides a unified approach to address a variety
of mathematical problems arising in various disciplines. In 2012, Censor et
al. \cite{Censor 2012} proposed the theory of split variational inequality
problems (SVIP) whereas Moudafi \cite{Moudafi 2011} generalized the concept
of SVIP to that of split monotone variational inclusions (SMVIP). The split
equilibrium problems is a special case of SMVIP. The SMVIP have already been
studied and successfully employed as a model in intensity-modulated
radiation therapy treatment planning, see \cite{Censor 2005,Censor 2006}.
Moreover, this formalism is also at the core of modeling of many inverse
problems arising for phase retrieval and other real-world problems; for
instance, in sensor networks in computerized tomography and data
compression; see, for example, \cite{Byrne 2002,Combettes 1996}. Some
methods have been proposed and analyzed to solve SEP together with the fixed
point problem in Hilbert spaces, see, for example \cite{Kazmi
2013,MAAKhan-LNA,MAAKhan-SB,Suantai FPTA 2016} and the references cited
therein.\newline

In 2013, Chang et al. \cite{Chang-2013} studied the split feasibility
problem for a total asymptotically strict pseudo contraction in infinitely
dimensional Hilbert spaces. In 2015, Ma and Wang \cite{Ma-Wang} established
strong convergence results for the split common fixed point problem of total
asymptotically strict pseudo contractions in Hilbert spaces. Quite recently,
some methods have been proposed and analyzed in \cite{MAAKhan-LNA,MAAKhan-SB}
for the split equilibrium problem. Inspired and motivated by the above
mentioned results and the ongoing research in this direction, we aim to
employ a hybrid shrinking projection algorithm to find a common element in
the set of solutions of a finite family of split equilibrium problems and
the set of common fixed points of a finite family of total asymptotically
strict pseudo contractions in Hilbert spaces. Our results can be viewed as a
generalization and improvement of various existing results in the current
literature.\\[1mm]

\section{Preliminaries}

This section is devoted to recall some definitions and results required in
the sequel.

Let $C$ be a nonempty closed convex subset of a Hilbert space $H_{1}.$ For
each $x\in H_{1}$, there exists a unique nearest point of $C,$ denoted by $%
P_{C}x,$ such that 
\begin{equation*}
\left\Vert x-P_{C}x\right\Vert \leq \left\Vert x-y\right\Vert \text{ for all 
}y\in C.
\end{equation*}%
Such a mapping $P_{C}:H_{1}\rightarrow C$ is known as a metric projection or
a nearest point projection of $H_{1}$ onto $C.$ Moreover, $P_{C}$ satisfies
nonexpansiveness in a Hilbert space and $\left\langle
x-P_{C}x,P_{C}x-y\right\rangle \geq 0$ for all $x,y\in C.$ It is remarked
that $P_{C}$ is firmly nonexpansive mapping from $H_{1}$ onto $C,$ that is,%
\begin{equation*}
\left\Vert P_{C}x-P_{C}y\right\Vert ^{2}\leq \left\langle
x-y,P_{C}x-P_{C}y~\right\rangle ,\text{ for all }x,y\in C.
\end{equation*}

Recall that a nonlinear mapping $A:C\rightarrow H_{1}$ is $\lambda $-inverse
strongly monotone if it satisfies%
\begin{equation*}
\left\langle x-y~,Ax-Ay\right\rangle \geq \lambda \left\Vert
Ax-Ay\right\Vert ^{2}.
\end{equation*}%
Note that, if $A:=I-T$ is a $\lambda $-inverse strongly monotone mapping,
then:\newline
(i): $A$ is a $\left( \frac{1}{\lambda }\right) $-Lipschitz continuous
mapping;\newline
(ii): if $T$ is a nonexpansive mapping, then $A$ is a $\left( \frac{1}{2}%
\right) $-inverse strongly monotone mapping;\newline
(iii): if $\eta \in (0,2\lambda ],$ then $I-\eta A$ is a nonexpansive
mapping.\newline

The following lemma collects some well-known equations in the context of a
real Hilbert space.\newline
\textbf{Lemma 2.1.} Let $H_{1}$ be a real Hilbert space, then:\newline
(i): $\left\Vert x-y\right\Vert ^{2}=\left\Vert x\right\Vert ^{2}-\left\Vert
y\right\Vert ^{2}-2\left\langle x-y,y\right\rangle ,$ for all $x,y\in H_{1};$%
\newline
(ii): $\left\Vert x+y\right\Vert ^{2}\leq \left\Vert x\right\Vert
^{2}+2\left\langle x-y,y\right\rangle ,$ for all $x,y\in H_{1};$\newline
(iii): $\left\Vert \alpha x+(1-\alpha )y\right\Vert ^{2}=\alpha \left\Vert
x\right\Vert ^{2}+(1-\alpha )\left\Vert y\right\Vert ^{2}-\alpha (1-\alpha
)\left\Vert x-y\right\Vert ^{2}$ for all $x,y\in H_{1}$ and $\alpha \in
\lbrack 0,1]$.\newline
\textbf{Lemma 2.2 }\cite{Moudafi 2010}\textbf{.} Let $T:C\rightarrow C$ be a 
$(k,\{\lambda _{n}\},\{\mu _{n}\},\xi )$-total asymptotically strictly
pseudo contraction. If $F(T)\neq \emptyset ,$ then for each $p\in F(T)$ and
for each $x\in C,$ the following equivalent inequalities hold: 
\begin{equation}
\left\Vert T^{n}x-p\right\Vert ^{2}\leq \left\Vert x-p\right\Vert
^{2}+k\left\Vert x-T^{n}x\right\Vert ^{2}+\lambda _{n}\xi (\Vert x-p\Vert
)+\mu _{n},\newline
\label{2.1}
\end{equation}%
\begin{equation}
\left\langle x-T^{n}x,x-p\right\rangle \geq \frac{1-k}{2}\left\Vert
x-T^{n}x\right\Vert ^{2}-\frac{\lambda _{n}}{2}\xi (\Vert x-p\Vert )-\frac{%
\mu _{n}}{2}\newline
,  \label{2.2}
\end{equation}%
\begin{equation}
\left\langle x-T^{n}x,p-T^{n}x\right\rangle \leq \frac{1+k}{2}\left\Vert
x-T^{n}x\right\Vert ^{2}+\frac{\lambda _{n}}{2}\xi (\Vert x-p\Vert )+\frac{%
\mu _{n}}{2}\newline
.  \label{2.3}
\end{equation}%
\newline
\textbf{Lemma 2.3 }\cite{Moudafi 2010}\textbf{.} Let $C$ be a nonempty
subset of a real Hilbert space $H_{1}$ and let $S:C\rightarrow C$ be a
uniformly $\Theta $-Lipschitzian and $(k,\{\lambda _{n}\},\{\mu _{n}\},\xi )$%
-total asymptotically strictly pseudo contraction, then $S$ is demiclosed at
origin. That is, if for any sequence $\{x_{n}\}$ in $C$ with $%
x_{n}\rightharpoonup x$ and $\left\Vert x_{n}-Sx_{n}\right\Vert \rightarrow
0,$ we have $x=Sx.$\newline
\textbf{Condition 2.4 }\cite{Blum,Combettes}\textbf{.} Let $f:C\times
C\rightarrow \mathbb{R}$ be a bifunction satisfying the following conditions:%
\newline
1. $f(x,x)=0\ $for\ all$\ x\in C;$\newline
2. $f$ is monotone, that is, $f(x,y)+f(y,x)\leq 0$ for\ all$\ x,y\in C;$%
\newline
3. $f$ is upper hemicontinuous, that is, for each $x,y,z\in C,$%
\begin{equation*}
\lim_{t\rightarrow 0}f(tz+(1-t)x,y)\leq f(x,y);
\end{equation*}%
4. for each $x\in C,$ the function $y\mapsto f(x,y)$ is convex and lower
semi-continuous.\newline
\textbf{Lemma 2.5 }\cite{Combettes}. Let $C$\ be a closed convex subset of a
real Hilbert space $H_{1}$\ and let $f:C\times C\rightarrow \mathbb{R}$ be a
bifunction satisfying Lemma 2.4. For $r>0$\ and $x\in H_{1},$\ there exists $%
z\in C$\ such that 
\begin{equation*}
F(z,y)+\frac{1}{r}\langle y-z,z-x\rangle \geq 0,\ \text{for all }y\in C.
\end{equation*}%
Moreover, define a mapping $T_{r}^{F}:H_{1}\rightarrow C$\ by 
\begin{equation*}
T_{r}^{f}(x)=\left\{ z\in C:f(z,y)+\frac{1}{r}\langle y-z,z-x\rangle \geq
0,\ \ \text{for all }y\in C\right\} ,
\end{equation*}%
for all $x\in H_{1}$. Then, the following hold:\newline
(i) $T_{r}^{f}$ is single-valued;\newline
(ii) $T_{r}^{f}$ is firmly nonexpansive, i.e., for every $x,y\in H,$%
\begin{equation*}
\left\Vert T_{r}^{f}x-T_{r}^{f}y\right\Vert ^{2}\leq \left\langle
T_{r}^{f}x-T_{r}^{f}y,x-y\right\rangle
\end{equation*}%
(ii) $F(T_{r}^{f})=EP(f);$\newline
(iv) $EP(f)$ is closed and convex.\newline

It is remarked that if $g:Q\times Q\rightarrow \mathbb{R}$ is a bifunction
satisfying Lemma 2.4, then for $s>0$\ and $w\in H_{2}$ we can define a
mapping: 
\begin{equation*}
T_{s}^{g}(w)=\left\{ d\in C:g(d,e)+\frac{1}{s}\langle e-d,d-w\rangle \geq
0,\ \ \text{for all }e\in Q\right\} ,
\end{equation*}%
which is, nonempty, single-valued and firmly nonexpansive. Moreover, $EP(g)$
is closed and convex, and $F(T_{s}^{g})=EP(g)$.

\section{Main results}

We now prove our main result of this section.\newline
\textbf{Theorem 3.1.} Let $H_{1}$ and $H_{2}$ be two real Hilbert spaces and
let $C\subseteq H_{1}$ and $Q\subseteq H_{2}$ be nonempty closed convex
subsets of Hilbert spaces $H_{1}$ and $H_{2}$, respectively. Let $%
f_{i}:C\times C\rightarrow \mathbb{R}$ and $g_{i}:Q\times Q\rightarrow 
\mathbb{R}$ be two finite families of bifunctions satisfying Condition 2.4
such that $g_{i}$ be upper semicontinuous for each $i\in \{1,2,3,\cdots ,N\}$%
. Let $S_{i}:C\rightarrow C$ be a finite family of uniformly $\Theta $%
-Lipschitzian and continuous\textit{\ }total asymptotically strict pseudo
contractions and let $A_{i}:H_{1}\rightarrow H_{2}$ be a finite family of
bounded linear operators for each $i\in \{1,2,3,\cdots ,N\}$. Suppose that $%
\mathbb{F}:=\left[ \bigcap_{i=1}^{N}F(S_{i})\right] \cap \Omega \neq
\emptyset $, where $\Omega =\left\{ z\in C:z\in \bigcap_{i=1}^{N}EP(f_{i})%
\text{ and }A_{i}z\in EP(g_{i})\text{ for }1\leq i\leq N\right\} .$ Let $%
\{x_{n}\}$\ be a sequence generated by:%
\begin{equation}
\begin{array}{l}
x_{1}\in C_{1}=C, \\ 
u_{n}=T_{r_{n}}^{f_{n}}\left( x_{n}-\gamma A_{n(\func{mod}N)}^{\ast }\left(
I-T_{s_{n}}^{g_{n}}\right) A_{n(\func{mod}N)}x_{n}\right) , \\ 
y_{n}=\alpha _{n}u_{n}+\left( 1-\alpha _{n}\right) S_{n(\func{mod}%
N)}^{n}u_{n}, \\ 
C_{n+1}=\left\{ {z\in H}_{1}{:}\left\Vert {y_{n}-z}\right\Vert {^{2}\leq
\left\Vert {x_{n}-z}\right\Vert ^{2}+\theta _{n}}\right\} , \\ 
x_{n+1}=P_{C_{n+1}}x_{1},\ \ n\geq 1,%
\end{array}
\label{3.1}
\end{equation}%
where $\theta _{n}=(1-\alpha _{n})\left\{ \lambda _{n}\xi
_{n}(M_{n})+\lambda _{n}M_{n}^{\ast }D_{n}+\mu _{n}\right\} $ with $%
D_{n}=\sup \left\{ \left\Vert x_{n}-p\right\Vert :p\ \in \mathbb{F}\right\} $%
. Let $\{r_{n}\},\{s_{n}\}$\ be two positive real sequences and let $%
\{\alpha _{n}\}$ be in $(0,1).$\ Assume that if the following set of
conditions holds:\newline
(C1): $0\leq k<a\leq \alpha _{n}\leq b<1$ and $\gamma \in \left( 0,\frac{1}{L%
}\right) $ where $L=\max \left\{ L_{1},L_{2},\cdots ,L_{N}\right\} $ and $%
L_{i}$ is the spectral radius of the operator $A_{i}^{\ast }A_{i}$ and $%
A_{i}^{\ast }$ is the adjoint of $A_{i}$ for each $i\in \{1,2,3,\cdots ,N\};$%
\newline
(C2): $\liminf\limits_{n\rightarrow \infty }r_{n}>0$ and $%
\liminf\limits_{n\rightarrow \infty }s_{n}>0;$\newline
(C3): $\sum\limits_{n=1}^{\infty }\lambda _{n}<\infty $ \textit{and}$\
\sum\limits_{n=1}^{\infty }\mu _{n}<\infty ;$\newline
(C4): \textit{there exist constants }$M_{i},\ M_{i}^{\ast }>0$ \textit{such
that }$\xi _{i}\left( \lambda _{i}\right) \leq M_{i}^{\ast }\lambda _{i}$%
\textit{\ for all }$\lambda _{i}\geq M_{i},i=1,2,3,\cdots ,N,$ then the
sequence $\{x_{n}\}$\ generated by (3.1) converges strongly to $P_{\mathbb{F}%
}x_{1}.$ \newline
\textbf{Proof. }For the sake of simplicity, we define $A_{n}=A_{n(\func{mod}%
N)}$ and $S_{n}=S_{n(\func{mod}N)}$ for all $n\geq 1.$ We start our proof to
establish that the sequence $\{x_{n}\}$\ defined in (3.1) is well defined.
In order to prove this assertion, we first show by mathematical induction
that $\mathbb{F}\subset C_{n}$ for all $n\geq 1.$ Obviously, $\mathbb{F}%
\subset C_{1}=C$. Now, assume that $\mathbb{F}\subset C_{i}$ for some $i\geq
1.$ Then it follows from (3.1) that%
\begin{eqnarray}
\left\Vert u_{i}-p\right\Vert ^{2} &=&\left\Vert T_{r_{i}}^{f_{i}}\left(
x_{i}-\gamma A_{i}^{\ast }\left( I-T_{s_{i}}^{g_{i}}\right)
A_{i}x_{i}\right) -T_{r_{i}}^{f_{i}}p\right\Vert ^{2}  \notag \\
&\leq &\left\Vert x_{i}-\gamma A_{i}^{\ast }\left(
I-T_{s_{i}}^{g_{i}}\right) A_{i}x_{i}-p\right\Vert ^{2}  \notag \\
&\leq &\left\Vert x_{i}-p\right\Vert ^{2}+\gamma ^{2}\left\Vert A_{i}^{\ast
}\left( I-T_{s_{i}}^{g_{i}}\right) A_{i}x_{i}\right\Vert ^{2}+2\gamma
\left\langle p-x_{i},A_{i}^{\ast }\left( I-T_{s_{i}}^{g_{i}}\right)
A_{i}x_{i}\right\rangle  \notag \\
&\leq &\left\Vert x_{i}-p\right\Vert ^{2}+\gamma ^{2}\left\langle
A_{i}x_{i}-T_{s_{i}}^{g_{i}}A_{i}x_{i},A_{i}A_{i}^{\ast }\left(
I-T_{s_{i}}^{g_{i}}\right) A_{i}x_{i}\right\rangle  \notag \\
&&+2\gamma \left\langle p-x_{i},A_{i}^{\ast }\left(
I-T_{s_{i}}^{g_{i}}\right) A_{i}x_{i}\right\rangle  \notag \\
&\leq &\left\Vert x_{i}-p\right\Vert ^{2}+L\gamma ^{2}\left\langle
A_{i}x_{i}-T_{s_{i}}^{g_{i}}A_{i}x_{i},A_{i}x_{i}-T_{s_{i}}^{g_{i}}A_{i}x_{i}\right\rangle
\notag \\
&&+2\gamma \left\langle p-x_{i},A_{i}^{\ast }\left(
I-T_{s_{i}}^{g_{i}}\right) A_{i}x_{i}\right\rangle  \notag \\
&=&\left\Vert x_{i}-p\right\Vert ^{2}+L\gamma ^{2}\left\Vert
A_{i}x_{i}-T_{s_{i}}^{g_{i}}A_{i}x_{i}\right\Vert ^{2}+2\gamma \left\langle
p-x_{i},A_{i}^{\ast }\left( I-T_{s_{i}}^{g_{i}}\right)
A_{i}x_{i}\right\rangle .  \label{3.2}
\end{eqnarray}%
Now letting $\Lambda =2\gamma \left\langle p-x_{i},A_{i}^{\ast }\left(
I-T_{s_{i}}^{g_{i}}\right) A_{i}x_{n}\right\rangle ,$ we have%
\begin{eqnarray*}
\Lambda &=&2\gamma \left\langle p-x_{i},A_{i}^{\ast }\left(
I-T_{s_{i}}^{g_{i}}\right) A_{i}x_{i}\right\rangle \\
&=&2\gamma \left\langle A_{i}\left( p-x_{i}\right)
,A_{i}x_{i}-T_{s_{i}}^{g_{i}}A_{i}x_{i}\right\rangle \\
&=&2\gamma \left\langle A_{i}\left( p-x_{i}\right) +\left(
A_{i}x_{i}-T_{s_{i}}^{g_{i}}A_{i}x_{i}\right) -\left(
A_{i}x_{i}-T_{s_{i}}^{g_{i}}A_{i}x_{i}\right)
,A_{i}x_{i}-T_{s_{i}}^{g_{i}}A_{i}x_{i}\right\rangle \\
&=&2\gamma \left\{ \left\langle
A_{i}p-T_{s_{i}}^{g_{i}}A_{i}x_{i},A_{i}x_{i}-T_{s_{i}}^{g_{i}}A_{i}x_{i}%
\right\rangle -\left\Vert A_{i}x_{i}-T_{s_{i}}^{g_{i}}A_{i}x_{i}\right\Vert
^{2}\right\} \\
&\leq &2\gamma \left\{ \frac{1}{2}\left\Vert
A_{i}x_{i}-T_{s_{i}}^{g_{i}}A_{i}x_{i}\right\Vert ^{2}-\left\Vert
A_{i}x_{i}-T_{s_{i}}^{g_{i}}A_{i}x_{i}\right\Vert ^{2}\right\} \\
&=&-\gamma \left\Vert A_{i}x_{i}-T_{s_{i}}^{g_{i}}A_{i}x_{i}\right\Vert ^{2}.
\end{eqnarray*}%
Using the above simplification of $\Lambda $ in (3.2), we get%
\begin{equation}
\left\Vert u_{i}-p\right\Vert ^{2}\leq \left\Vert x_{i}-p\right\Vert
^{2}+\gamma \left( L\gamma -1\right) \left\Vert
A_{i}x_{i}-T_{s_{n,i}}^{g_{i}}A_{i}x_{i}\right\Vert ^{2}.  \label{3.3}
\end{equation}%
Since $\gamma \in \left( 0,\frac{1}{L}\right) $ by condition (C1), the above
estimate then yields%
\begin{equation}
\left\Vert u_{i}-p\right\Vert ^{2}\leq \left\Vert x_{i}-p\right\Vert ^{2}.
\label{3.4}
\end{equation}%
Making use of (3.4), we have the following estimate:%
\begin{eqnarray}
\left\Vert y_{i}-p\right\Vert ^{2} &=&\left\Vert \alpha _{i}u_{i}+\left(
1-\alpha _{i}\right) S_{i}^{i}u_{i}-p\right\Vert ^{2}  \notag \\
&=&\alpha _{i}\left\Vert u_{i}-p\right\Vert ^{2}+\left( 1-\alpha _{i}\right)
\left\Vert S_{i}^{i}u_{i}-p\right\Vert ^{2}-\alpha _{i}\left( 1-\alpha
_{i}\right) \left\Vert u_{i}-S_{i}^{i}u_{i}\right\Vert ^{2}  \notag \\
&\leq &\alpha _{i}\left\Vert u_{i}-p\right\Vert ^{2}+\left( 1-\alpha
_{i}\right) \left\{ \left\Vert u_{i}-p\right\Vert ^{2}+k\left\Vert
u_{i}-S_{i}^{i}u_{i}\right\Vert ^{2}+\lambda _{i}\xi _{i}(\left\Vert
u_{i}-p\right\Vert )+\mu _{i}\right\}  \notag \\
&&-\alpha _{i}\left( 1-\alpha _{i}\right) \left\Vert
u_{i}-S_{i}^{i}u_{i}\right\Vert ^{2}  \notag \\
&\leq &\left\Vert u_{i}-p\right\Vert ^{2}+\left( 1-\alpha _{i}\right)
\left\{ k\left\Vert u_{i}-S_{i}^{i}u_{i}\right\Vert ^{2}+\lambda _{i}\xi
_{i}\left( M_{i}\right) +\lambda _{i}M_{i}^{\ast }\left\Vert
u_{i}-p\right\Vert ^{2})+\mu _{i}\right\}  \notag \\
&&-\alpha _{i}\left( 1-\alpha _{i}\right) \left\Vert
u_{i}-S_{i}^{i}u_{i}\right\Vert ^{2}  \notag \\
&\leq &\left\Vert x_{i}-p\right\Vert ^{2}-\left( 1-\alpha _{i}\right) \left(
\alpha _{i}-k\right) \left\Vert u_{i}-S_{i}^{i}u_{i}\right\Vert ^{2}  \notag
\\
&&+\left( 1-\alpha _{i}\right) \left\{ \lambda _{i}\xi _{i}\left(
M_{i}\right) +\lambda _{i}M_{i}^{\ast }\left\Vert x_{i}-p\right\Vert
^{2}+\mu _{i}\right\} .  \label{3.5}
\end{eqnarray}%
Since $\alpha _{i}-k\geq 0~$by condition (C1), so (3.5) implies that%
\begin{equation}
\left\Vert y_{i}-p\right\Vert ^{2}\leq \left\Vert x_{i}-p\right\Vert
^{2}+\theta _{i},  \label{3.6}
\end{equation}%
\newline
where $\theta _{i}=\left( 1-\alpha _{i}\right) \left\{ \lambda _{i}\xi
_{i}(M_{i})+\lambda _{i}M_{i}^{\ast }D_{i}+\mu _{i}\right\} $ and $%
D_{i}=\sup \left\{ \left\Vert x_{i}-p\right\Vert ^{2}:p\ \in \mathbb{F}%
\right\} .$ It now follows from the estimate (3.6) that $p\in C_{i+1}.$
Hence, $\mathbb{F}\subset C_{n}$ for all $n\geq 1.$ Next, we show that the
set $C_{n}$ is closed and convex for all $n\geq 1.$ Since%
\begin{equation*}
\left\{ {z\in H_{1}:}\left\Vert {y_{n}-z}\right\Vert {^{2}\leq \left\Vert {%
x_{n}-z}\right\Vert ^{2}}\right\} =\left\{ {z\in H_{1}:}\left\Vert {y_{n}}%
\right\Vert {^{2}-\left\Vert {x_{n}}\right\Vert {^{2}}\leq 2}\left\langle {%
y_{n}}-x_{n},z\right\rangle {+\theta _{n}}\right\} ,
\end{equation*}%
it is closed and convex; hence the sequence $\left\{ x_{n}\right\} $ defined
in (3.1) is well-defined. Next, from $x_{n}=P_{C_{n}}x_{1},$ we get%
\begin{equation*}
0\leq \left\langle x_{n}-x_{1},x^{\ast }-x_{n}\right\rangle ,
\end{equation*}%
for all $x^{\ast }\in C_{n}.$ Since $\mathbb{F\neq \emptyset },$ then for
any $p\in \mathbb{F\,}$, we get%
\begin{eqnarray*}
0 &\leq &\left\langle x_{n}-x_{1},p-x_{n}\right\rangle \\
&=&\left\langle x_{n}-x_{1},p+x_{1}-x_{1}-x_{n}\right\rangle \\
&=&\left\langle x_{n}-x_{1},x_{1}-x_{n}\right\rangle +\left\langle
x_{n}-x_{1},p-x_{1}\right\rangle \\
&=&-\left\Vert x_{n}-x_{1}\right\Vert ^{2}+\left\Vert x_{n}-x_{1}\right\Vert
\left\Vert p-x_{1}\right\Vert .
\end{eqnarray*}%
That is,%
\begin{equation*}
\left\Vert x_{n}-x_{1}\right\Vert \leq \left\Vert p-x_{1}\right\Vert ,\text{
for all }p\in \mathbb{F}\text{ and }n\geq 1.
\end{equation*}%
Hence, the sequence $\{x_{n}\}$ is bounded, so are $\{u_{n}\}$ and $%
\{y_{n}\}.$ Moreover, from $x_{n}=P_{C_{n}}x_{1}$ and $%
x_{n+1}=P_{C_{n+1}}x_{1}\in C_{n+1}\subset C_{n},$ we have%
\begin{equation*}
0\leq \left\langle x_{n}-x_{1},x_{n+1}-x_{n}\right\rangle .
\end{equation*}%
Similarly, we get the following relation:%
\begin{equation*}
\left\Vert x_{n}-x_{1}\right\Vert \leq \left\Vert x_{n+1}-x_{1}\right\Vert ,%
\text{ for all }n\geq 1.
\end{equation*}%
From the above assertions, we conclude that the sequence $\{\left\Vert
x_{n}-x_{1}\right\Vert \}$ is bounded and nondecreasing, therefore, we have%
\begin{equation}
\lim_{n\rightarrow \infty }\left\Vert x_{n}-x_{1}\right\Vert \text{ exists.}
\label{3.7}
\end{equation}%
Further observe that%
\begin{eqnarray*}
\left\Vert x_{n+1}-x_{n}\right\Vert ^{2} &=&\left\Vert
x_{n+1}-x_{1}+x_{1}-x_{n}\right\Vert ^{2} \\
&=&\left\Vert x_{n+1}-x_{1}\right\Vert ^{2}+\left\Vert
x_{n}-x_{1}\right\Vert ^{2}-2\left\langle
x_{n}-x_{1},x_{n+1}-x_{1}\right\rangle \\
&=&\left\Vert x_{n+1}-x_{1}\right\Vert ^{2}+\left\Vert
x_{n}-x_{1}\right\Vert ^{2}-2\left\langle
x_{n}-x_{1},x_{n+1}-x_{n}+x_{n}-x_{1}\right\rangle \\
&=&\left\Vert x_{n+1}-x_{1}\right\Vert ^{2}-\left\Vert
x_{n}-x_{1}\right\Vert ^{2}-2\left\langle
x_{n}-x_{1},x_{n+1}-x_{n}\right\rangle \\
&\leq &\left\Vert x_{n+1}-x_{1}\right\Vert ^{2}-\left\Vert
x_{n}-x_{1}\right\Vert ^{2}.
\end{eqnarray*}%
From (3.7), we obtain that $\left\Vert x_{n+1}-x_{1}\right\Vert
^{2}-\left\Vert x_{n}-x_{1}\right\Vert ^{2}\rightarrow 0$ as $n\rightarrow
\infty ,$ hence the sequence $\{\left\Vert x_{n}-x_{1}\right\Vert \}$ is
Cauchy. That is%
\begin{equation}
\lim_{n\rightarrow \infty }\left\Vert x_{n+1}-x_{n}\right\Vert =0.
\label{3.8}
\end{equation}%
Since $x_{n+1}\in C_{n+1},$ which implies that $\left\Vert
y_{n}-x_{n+1}\right\Vert \leq \left\Vert x_{n}-x_{n+1}\right\Vert +\theta
_{n,i}.$ From (3.8), we conclude that%
\begin{equation}
\lim_{n\rightarrow \infty }\left\Vert y_{n}-x_{n+1}\right\Vert =0,\text{ for
all }n\geq 1.  \label{3.9}
\end{equation}%
Now using (3.8), (3.9) and the following triangular inequality, we get%
\begin{equation}
\left\Vert y_{n}-x_{n}\right\Vert \leq \left\Vert y_{n}-x_{n+1}\right\Vert
+\left\Vert x_{n+1}-x_{n}\right\Vert \rightarrow 0  \label{3.10}
\end{equation}%
as $n\rightarrow \infty .$\newline
Consider from (3.1), (3.3) and (3.6), we get%
\begin{eqnarray*}
\gamma \left( 1-\gamma L\right) \left\Vert
A_{n}x_{n}-T_{s_{n}}^{g_{n}}A_{n}x_{n}\right\Vert ^{2} &\leq &\left\Vert
x_{n}-p\right\Vert ^{2}-\left\Vert u_{n}-p\right\Vert ^{2} \\
&\leq &\left\Vert x_{n}-p\right\Vert ^{2}-\left\Vert y_{n}-p\right\Vert
^{2}+\theta _{n} \\
&\leq &\left( \left\Vert x_{n}-p\right\Vert +\left\Vert y_{n}-p\right\Vert
\right) \left\Vert x_{n}-y_{n}\right\Vert +\theta _{n}.
\end{eqnarray*}%
Utilizing the fact that $\gamma \left( 1-\gamma L\right) >0$ and the
estimate (3.10), we have%
\begin{equation}
\lim_{n\rightarrow \infty }\left\Vert
A_{n}x_{n}-T_{s_{n}}^{g_{n}}A_{n}x_{n}\right\Vert ^{2}=0\text{ for all }%
n\geq 1.  \label{3.11}
\end{equation}%
For any $p\in \mathbb{F}$ and firm nonexpansiveness of $T_{r_{n}}^{f_{n}},$
we have%
\begin{eqnarray}
\left\Vert u_{n}-p\right\Vert ^{2} &=&\left\Vert T_{r_{n}}^{f_{n}}\left(
x_{n}-\gamma A_{n}^{\ast }\left( I-T_{s_{n}}^{g_{n}}\right)
A_{n}x_{n}\right) -T_{r_{n}}^{f_{n}}p\right\Vert ^{2}  \notag \\
&\leq &\left\langle u_{n}-p,x_{n}-\gamma A_{n}^{\ast }\left(
I-T_{s_{n}}^{g_{n}}\right) A_{n}x_{n}-p\right\rangle  \notag \\
&=&\frac{1}{2}\{\left\Vert u_{n}-p\right\Vert ^{2}+\left\Vert x_{n}-\gamma
A_{n}^{\ast }\left( I-T_{s_{n}}^{g_{n}}\right) A_{n}x_{n}-p\right\Vert ^{2} 
\notag \\
&&-\left\Vert u_{n}-x_{n}-\gamma A_{n}^{\ast }\left(
I-T_{s_{n}}^{g_{n}}\right) A_{n}x_{n}\right\Vert ^{2}\}  \notag \\
&\leq &\frac{1}{2}\left\{ \left\Vert u_{n}-p\right\Vert ^{2}+\left\Vert
x_{n}-p\right\Vert ^{2}-\left\Vert u_{n}-x_{n}-\gamma A_{n}^{\ast }\left(
I-T_{s_{n}}^{g_{n}}\right) A_{n}x_{n}\right\Vert ^{2}\right\}  \notag \\
&=&\frac{1}{2}\{\left\Vert u_{n}-p\right\Vert ^{2}+\left\Vert
x_{n}-p\right\Vert ^{2}-(\left\Vert u_{n}-x_{n}\right\Vert ^{2}+\gamma
^{2}\left\Vert A_{n}^{\ast }\left( I-T_{s_{n}}^{g_{n}}\right)
A_{n}x_{n}\right\Vert ^{2}  \notag \\
&&-2\gamma \left\langle u_{n}-x_{n},A_{n}^{\ast }\left(
I-T_{s_{n}}^{g_{n}}\right) A_{n}x_{n}\right\rangle )\}  \notag \\
&\leq &\left\Vert x_{n}-p\right\Vert ^{2}-\left\Vert u_{n}-x_{n}\right\Vert
^{2}+2\gamma \left\Vert A\left( u_{n}-x_{n}\right) \right\Vert \left\Vert
A_{n}x_{n}-T_{s_{n}}^{g_{n}}A_{n}x_{n}\right\Vert .  \label{3.12}
\end{eqnarray}%
Using the following estimate:%
\begin{equation*}
\left\Vert y_{n}-p\right\Vert ^{2}\leq \alpha _{n}\left\Vert
x_{n}-p\right\Vert ^{2}+\left( 1-\alpha _{n}\right) \left\Vert
u_{n}-p\right\Vert ^{2}+\theta _{n}.
\end{equation*}%
in (3.12) and re-arranging the terms, we get%
\begin{eqnarray*}
\left( 1-\alpha _{n}\right) \left\Vert u_{n}-x_{n}\right\Vert ^{2} &\leq
&\left( \left\Vert x_{n}-p\right\Vert +\left\Vert y_{n}-p\right\Vert \right)
\left\Vert x_{n}-y_{n}\right\Vert \\
&&+2\gamma \left\Vert A\left( u_{n}-x_{n}\right) \right\Vert \left\Vert
A_{n}x_{n}-T_{s_{n}}^{g_{n}}A_{n}x_{n}\right\Vert +\theta _{n}.
\end{eqnarray*}%
Letting $n\rightarrow \infty $ and utilizing (3.10) and (3.11), we have%
\begin{equation}
\lim_{n\rightarrow \infty }\left\Vert u_{n}-x_{n}\right\Vert =0\text{ for
all }n\geq 1.  \label{3.13}
\end{equation}%
\newline
Moreover, from (3.10) and (3.13), we obtain%
\begin{equation}
\left\Vert y_{n}-u_{n}\right\Vert \leq \left\Vert y_{n}-x_{n}\right\Vert
+\left\Vert x_{n}-u_{n}\right\Vert \rightarrow 0  \label{3.14}
\end{equation}%
when $n\rightarrow \infty .$ Observe that $\left\Vert y_{n}-u_{n}\right\Vert
=\left( 1-\alpha _{n}\right) \left\Vert S_{n}^{n}u_{n}-u_{n}\right\Vert .$
Then it follows from condition (C1) and (3.14) that%
\begin{equation}
\lim_{n\rightarrow \infty }\left\Vert S_{n}^{n}u_{n}-u_{n}\right\Vert =0%
\text{ for all }n\geq 1.  \label{3.15}
\end{equation}%
Since%
\begin{equation*}
\left\Vert S_{n}^{n}u_{n}-x_{n}\right\Vert \leq \left\Vert
S_{n}^{n}u_{n}-u_{n}\right\Vert +\left\Vert x_{n}-u_{n}\right\Vert .
\end{equation*}%
Therefore from (3.13) and (3.15), we get%
\begin{equation}
\lim_{n\rightarrow \infty }\left\Vert S_{n}^{n}u_{n}-x_{n}\right\Vert =0%
\text{ for all }n\geq 1.  \label{3.16}
\end{equation}%
On a similar reasoning, we also obtain%
\begin{equation*}
\lim_{n\rightarrow \infty }\left\Vert S_{n}^{n}u_{n}-y_{n}\right\Vert =0%
\text{ for all }n\geq 1.
\end{equation*}%
Observe that each $S_{n}$ is uniformly $\Theta $-Lipschitzian, therefore, we
have%
\begin{eqnarray*}
\left\Vert S_{n}^{n}x_{n}-x_{n}\right\Vert &\leq &\left\Vert
S_{n}^{n}x_{n}-S_{n}^{n}u_{n}\right\Vert +\left\Vert
S_{n}^{n}u_{n}-x_{n}\right\Vert \\
&\leq &\Theta \left\Vert x_{n}-u_{n}\right\Vert +\left\Vert
S_{n}^{n}u_{n}-x_{n}\right\Vert .
\end{eqnarray*}%
Now, using (3.13), (3.16) and the above estimate, we get%
\begin{equation}
\lim_{n\rightarrow \infty }\left\Vert S_{n}^{n}x_{n}-x_{n}\right\Vert =0%
\text{ for all }n\geq 1.  \label{3.17}
\end{equation}%
Moreover, utilizing the uniform continuity of $S_{n}$ and (3.17)$,$ the
following estimate:%
\begin{equation*}
\left\Vert x_{n}-S_{n}x_{n}\right\Vert \leq \left\Vert
x_{n}-S_{n}^{n}x_{n}\right\Vert +\left\Vert
S_{n}^{n}x_{n}-S_{n}x_{n}\right\Vert
\end{equation*}%
implies that%
\begin{equation}
\lim_{n\rightarrow \infty }\left\Vert x_{n}-S_{n}x_{n}\right\Vert =0\text{
for all }n\geq 1.  \label{3.18}
\end{equation}%
Similarly, we also have that%
\begin{equation}
\lim_{n\rightarrow \infty }\left\Vert u_{n}-S_{n}u_{n}\right\Vert =0\text{
for all }n\geq 1.  \label{3.19}
\end{equation}%
Now, we show that $\omega (x_{n})\subset \mathbb{F},$ where $\omega (x_{n})$
is the set of all weak $\omega $-limits of $\{x_{n}\}.$ Since $\{x_{n}\}$ is
bounded, therefore $\omega (x_{n})\neq \emptyset .$ Let $q\in \omega
(x_{n}), $ then there exists a subsequence $\{x_{Nn+i}\}$ of $\{x_{n}\}$
such that $x_{Nn+i}\rightharpoonup q.$ Using the fact that $S_{Nn+i}=S_{i}$
for all $n\geq 1$ and the demiclosed principle (Lemma 2.3) for each $S_{i}$,
we have that $x\in F(S_{i})\,$\ for each $1\leq i\leq N.$ Next, we show that 
$q\in \Omega $, i.e., $q\in \bigcap_{i=1}^{N}EP(f_{i})$ and $A_{i}q\in
EP(g_{i})$ for each $1\leq i\leq N.$ In order to show that $q\in
\bigcap_{i=1}^{N}EP(f_{i}),$ that is, $q\in EP(f_{i})$ for each $1\leq i\leq
N,$ we define subsequence $\left\{ n_{j}\right\} $ of index $\left\{
n\right\} $ such that $n_{j}=Nj+i$ for all $n\geq 1.$ As a consequence, we
can write $f_{n_{j}}=f_{i}$ for $1\leq i\leq N.$ From $%
u_{n_{j}}=T_{r_{n_{j}}}^{f_{i}}\left( I-\gamma A_{n_{j}}^{\ast }\left(
I-T_{s_{n_{j}}}^{g_{n_{j}}}\right) A_{n_{j}}\right) x_{n_{j}}$ for all $%
n\geq 1,$ we have%
\begin{equation*}
f_{i}(u_{n_{j}},y)+\frac{1}{r_{n_{j}}}\left\langle
y-u_{n_{j}},u_{n_{j}}-x_{n_{j}}-\gamma A_{n_{j}}^{\ast }\left(
I-T_{s_{n_{j}}}^{g_{n_{j}}}\right) A_{n_{j}}x_{n_{j}}\right\rangle \geq 0,\
\ \text{for all }y\in C.
\end{equation*}%
This implies that%
\begin{equation*}
f_{i}(u_{n_{j}},y)+\frac{1}{r_{n_{j}}}\left\langle
y-u_{n_{j}},u_{n_{j}}-x_{n_{j}}\right\rangle -\frac{1}{r_{n_{j}}}%
\left\langle y-u_{n_{j}},\gamma A_{n_{j}}^{\ast }\left(
I-T_{s_{n_{j}}}^{g_{n_{j}}}\right) A_{n_{j}}x_{n_{j}}\right\rangle \geq 0
\end{equation*}%
From (A2), we have%
\begin{equation*}
\frac{1}{r_{n_{j}}}\left\langle y-u_{n_{j}},u_{n_{j}}-x_{n_{j}}\right\rangle
-\frac{1}{r_{n_{j}}}\left\langle y-u_{n_{j}},\gamma A_{n_{j}}^{\ast }\left(
I-T_{s_{n_{j}}}^{g_{n_{j}}}\right) A_{n_{j}}x_{n_{j}}\right\rangle \geq
f_{i}(y,u_{n_{j}}),
\end{equation*}%
for all $y\in C.$ Since $\liminf_{j\rightarrow \infty }r_{n_{i}}>0$ (by
(C2)), therefore it follows from (3.11) and (3.13) that%
\begin{equation*}
f_{i}(y,q)\leq 0,\ \ \text{for all }y\in C\text{ and for }1\leq i\leq N.
\end{equation*}%
Let $y_{t}=ty+(1-t)q$ for some $0<t<1$ and $y\in C$. Since $q\in C,$ this
implies that $y_{t}\in C.$ Using (A1) and (A4) from Condition 2.4, the
following estimate: 
\begin{equation*}
0=f_{i}(y_{t},y_{t})\leq tf_{i}(y_{t},y)+(1-t)f_{i}(y_{t},q)\leq
tf_{i}(y_{t},y),
\end{equation*}%
implies that%
\begin{equation*}
f_{i}(y_{t},y)\geq 0,\text{ for }1\leq i\leq N.
\end{equation*}%
Letting $t\rightarrow 0,$ we have $f_{i}(q,y)\geq 0$ for all $y\in C.$ Thus, 
$q\in EP(f_{i})$ for $1\leq i\leq N.$ That is, $q\in
\bigcap_{i=1}^{N}EP(F_{i}).$ Reasoning as above, we show that $A_{i}q\in
EP(g_{i})$ for each $1\leq i\leq N.$ Since $x_{n_{l}}\longrightarrow q$ and $%
A_{n_{l}}$ is a bounded linear operator, therefore $A_{n_{l}}x_{n_{l}}%
\longrightarrow A_{n_{l}}q.$ Hence, it follows from (3.11) that%
\begin{equation*}
T_{s_{n_{l}}}^{g_{n_{l}}}A_{n_{l}}x_{n_{l}}\longrightarrow A_{n_{l}}q\text{
\ \ \ as \ \ }l\rightarrow \infty .
\end{equation*}%
Now, from Lemma 2.5, we have%
\begin{equation*}
g_{i}\left( T_{s_{n_{l}}}^{g_{n_{l}}}A_{n_{l}}x_{n_{l}},z\right) +\frac{1}{%
s_{n_{l}}}\left\langle
z-T_{s_{n_{l}}}^{g_{n_{l}}}A_{n_{l}}x_{n_{l}},T_{s_{n_{l}}}^{g_{n_{l}}}A_{n_{l}}x_{n_{l}}-A_{n_{l}}x_{n_{l}}\right\rangle \geq 0,\ \ 
\text{for all }z\in Q.
\end{equation*}%
Since $g_{i}$ is upper hemicontinuous in the first argument for each $1\leq
i\leq N$, therefore taking $\limsup $ on both sides of the above estimate as 
$l\rightarrow \infty $ and utilizing (C2) and (3.11), we get%
\begin{equation*}
g_{i}\left( A_{n_{l}}x,z\right) \geq 0,\ \ \text{for all }z\in Q\text{ and
for each }1\leq i\leq N.
\end{equation*}%
Hence $A_{i}q\in EP(g_{i})$ for each $1\leq i\leq N$ and consequently $q\in 
\mathbb{F}.$ It remains to show that $x_{n}\rightarrow q=P_{\mathbb{F}%
}x_{1}. $ Let $x=P_{\mathbb{F}}x_{1},$ then from $\left\Vert
x_{n}-x_{1}\right\Vert \leq \left\Vert x-x_{1}\right\Vert ,$ therefore, we
have%
\begin{eqnarray*}
\left\Vert x-x_{1}\right\Vert &\leq &\left\Vert q-x_{1}\right\Vert \\
&\leq &\liminf_{j\rightarrow \infty }\left\Vert x_{n_{j}}-x_{1}\right\Vert \\
&\leq &\limsup_{j\rightarrow \infty }\left\Vert x_{n_{j}}-x_{1}\right\Vert \\
&\leq &\left\Vert x-x_{1}\right\Vert .
\end{eqnarray*}%
This implies that%
\begin{equation*}
\lim_{j\rightarrow \infty }\left\Vert x_{n_{j}}-x_{1}\right\Vert =\left\Vert
q-x_{1}\right\Vert .
\end{equation*}%
Hence $x_{n_{j}}\rightarrow q=P_{\mathbb{F}}x_{1}.$ From the arbitrariness
of the subsequence $\left\{ x_{n_{j}}\right\} $ of $\left\{ x_{n}\right\} ,$
we conclude that $x_{n}\rightarrow x$ as $n\rightarrow \infty .$ It is easy
to see that $y_{n,i}\rightarrow x$ and $u_{n,i}\rightarrow x.$ This
completes the proof. \ \ \ \ \ \ \ \ \ \ \ \ \ \ \ \ \ \ \ \ \ \ \ \ \ \ \ \
\ \ \ \ \ \ \ \ \ \ \ \ \ \ \ \ \ \ \ \ \ \ \ \ \ \ \ \ \ \ \ \ \ \ \ \ \ \
\ \ \ \ \ \ \ \ \ \ \ \ \ \ \ \ \ \ \ \ \ \ \ \ \ \ \ \ \ \ $\square $%
\newline
\textbf{Corollary 3.2.} Let $H_{1}$ and $H_{2}$ be two real Hilbert spaces
and let $C\subseteq H_{1}$ and $Q\subseteq H_{2}$ be nonempty closed convex
subsets of Hilbert spaces $H_{1}$ and $H_{2}$, respectively. Let $%
f_{i}:C\times C\rightarrow \mathbb{R}$ and $g_{i}:Q\times Q\rightarrow 
\mathbb{R}$ be two finite families of bifunctions satisfying Condition 2.4
such that $g_{i}$ be upper semicontinuous for each $i\in \{1,2,3,\cdots ,N\}$%
. Let $S_{i}:C\rightarrow C$ be a finite family nonexpansive mappings and
let $A_{i}:H_{1}\rightarrow H_{2}$ be a finite family of bounded linear
operators for each $i\in \{1,2,3,\cdots ,N\}$. Suppose that $\mathbb{F}:=%
\left[ \bigcap_{i=1}^{N}F(S_{i})\right] \cap \Omega \neq \emptyset $, where $%
\Omega =\left\{ z\in C:z\in \bigcap_{i=1}^{N}EP(f_{i})\text{ and }A_{i}z\in
EP(g_{i})\text{ for }1\leq i\leq N\right\} .$ Let $\{x_{n}\}$\ be a sequence
generated by:%
\begin{equation}
\begin{array}{l}
x_{1}\in C_{1}=C, \\ 
u_{n}=T_{r_{n}}^{f_{n}}\left( x_{n}-\gamma A_{n(\func{mod}N)}^{\ast }\left(
I-T_{s_{n}}^{g_{n}}\right) A_{n(\func{mod}N)}x_{n}\right) , \\ 
y_{n}=\alpha _{n}u_{n}+\left( 1-\alpha _{n}\right) S_{n(\func{mod}%
N)}^{n}u_{n}, \\ 
C_{n+1}=\left\{ {z\in H}_{1}{:}\left\Vert {y_{n}-z}\right\Vert {^{2}\leq
\left\Vert {x_{n}-z}\right\Vert ^{2}+\theta _{n}}\right\} , \\ 
x_{n+1}=P_{C_{n+1}}x_{1},\ \ n\geq 1,%
\end{array}
\label{3.20}
\end{equation}%
where $\theta _{n}=(1-\alpha _{n})\left\{ \lambda _{n}\xi
_{n}(M_{n})+\lambda _{n}M_{n}^{\ast }D_{n}+\mu _{n}\right\} $ with $%
D_{n}=\sup \left\{ \left\Vert x_{n}-p\right\Vert :p\ \in \mathbb{F}\right\} $%
. Let $\{r_{n}\},\{s_{n}\}$\ be two positive real sequences and let $%
\{\alpha _{n}\}$ be in $(0,1).$\ Assume that if the following set of
conditions holds:\newline
(C1): $0\leq k<a\leq \alpha _{n}\leq b<1$ and $\gamma \in \left( 0,\frac{1}{L%
}\right) $ where $L=\max \left\{ L_{1},L_{2},\cdots ,L_{N}\right\} $ and $%
L_{i}$ is the spectral radius of the operator $A_{i}^{\ast }A_{i}$ and $%
A_{i}^{\ast }$ is the adjoint of $A_{i}$ for each $i\in \{1,2,3,\cdots ,N\};$%
\newline
(C2): $\liminf\limits_{n\rightarrow \infty }r_{n}>0$ and $%
\liminf\limits_{n\rightarrow \infty }s_{n}>0;$\newline
(C3): $\sum\limits_{n=1}^{\infty }\lambda _{n}<\infty $ \textit{and}$\
\sum\limits_{n=1}^{\infty }\mu _{n}<\infty ;$\newline
(C4): \textit{there exist constants }$M_{i},\ M_{i}^{\ast }>0$ \textit{such
that }$\xi _{i}\left( \lambda _{i}\right) \leq M_{i}^{\ast }\lambda _{i}$%
\textit{\ for all }$\lambda _{i}\geq M_{i},i=1,2,3,\cdots ,N,$ then the
sequence $\{x_{n}\}$\ generated by (3.20) converges strongly to $P_{\mathbb{F%
}}x_{1}.\medskip $\newline
\textbf{Theorem 3.3.}\ Let $H_{1}$ and $H_{2}$ be two real Hilbert spaces
and let $C\subseteq H_{1}$ and $Q\subseteq H_{2}$ be nonempty closed convex
subsets of Hilbert spaces $H_{1}$ and $H_{2}$, respectively. Let $%
f_{i}:C\times C\rightarrow \mathbb{R}$ and $g_{i}:Q\times Q\rightarrow 
\mathbb{R}$ be two finite families of bifunctions satisfying Condition 2.4
such that $g_{i}$ be upper semicontinuous for each $i\in \{1,2,3,\cdots ,N\}$%
. Let $S_{i}:C\rightarrow C$ be a finite family of uniformly $\Theta $%
-Lipschitzian and continuous\textit{\ }total asymptotically strict pseudo
contractions and let $A_{i}:H_{1}\rightarrow H_{2}$ be a finite family of
bounded linear operators for each $i\in \{1,2,3,\cdots ,N\}$. Suppose that $%
\mathbb{F}:=\left[ \bigcap_{i=1}^{N}F(S_{i})\right] \cap \Omega \neq
\emptyset $, where $\Omega =\left\{ z\in C:z\in \bigcap_{i=1}^{N}EP(f_{i})%
\text{ and }A_{i}z\in EP(g_{i})\text{ for }1\leq i\leq N\right\} .$ Let $%
\{x_{n}\}$\ be a sequence generated by:%
\begin{equation}
\begin{array}{l}
x_{1}\in C_{1}=C, \\ 
u_{n}=T_{r_{n}}^{f_{n}}\left( x_{n}-\gamma A_{n(\func{mod}N)}^{\ast }\left(
I-T_{s_{n}}^{g_{n}}\right) A_{n(\func{mod}N)}x_{n}\right) , \\ 
y_{n}=\alpha _{n}u_{n}+\left( 1-\alpha _{n}\right) S_{n(\func{mod}%
N)}^{n}u_{n}, \\ 
C_{n+1}=\left\{ {z\in H}_{1}{:}\left\Vert {y_{n}-z}\right\Vert {^{2}\leq
\left\Vert {x_{n}-z}\right\Vert ^{2}+\theta _{n}}\right\} , \\ 
x_{n+1}=P_{C_{n+1}}x_{1},\ \ n\geq 1,%
\end{array}
\label{3.21}
\end{equation}%
where $\theta _{n}=(1-\alpha _{n})\left\{ \lambda _{n}\xi
_{n}(M_{n})+\lambda _{n}M_{n}^{\ast }D_{n}+\mu _{n}\right\} $ with $%
D_{n}=\sup \left\{ \left\Vert x_{n}-p\right\Vert :p\ \in \mathbb{F}\right\} $%
. Let $\{r_{n}\},\{s_{n}\}$\ be two positive real sequences and let $%
\{\alpha _{n}\}$ be in $(0,1).$\ Assume that if the following set of
conditions holds:\newline
(C1): $0\leq k<a\leq \alpha _{n}\leq b<1$ and $\gamma \in \left( 0,\frac{1}{L%
}\right) $ where $L=\max \left\{ L_{1},L_{2},\cdots ,L_{N}\right\} $ and $%
L_{i}$ is the spectral radius of the operator $A_{i}^{\ast }A_{i}$ and $%
A_{i}^{\ast }$ is the adjoint of $A_{i}$ for each $i\in \{1,2,3,\cdots ,N\};$%
\newline
(C2): $\liminf\limits_{n\rightarrow \infty }r_{n}>0$ and $%
\liminf\limits_{n\rightarrow \infty }s_{n}>0;$\newline
(C3): $\sum\limits_{n=1}^{\infty }\lambda _{n}<\infty $ \textit{and}$\
\sum\limits_{n=1}^{\infty }\mu _{n}<\infty ;$\newline
(C4): \textit{there exist constants }$M_{i},\ M_{i}^{\ast }>0$ \textit{such
that }$\xi _{i}\left( \lambda _{i}\right) \leq M_{i}^{\ast }\lambda _{i}$%
\textit{\ for all }$\lambda _{i}\geq M_{i},i=1,2,3,\cdots ,N,$ then the
sequence $\{x_{n}\}$\ generated by (3.21) converges strongly to $P_{\mathbb{F%
}}x_{1}.$\newline
\textbf{Proof.} Set $H_{1}=H_{2},~C=Q$ and $A=I$(the identity mapping) then
the desired result then follows from Theorem 3.1 immediately.\ \ \ \ \ \ \ \
\ \ \ \ \ \ \ \ \ \ \ \ \ \ \ \ \ \ \ \ \ \ \ \ \ \ \ \ \ \ \ \ \ \ \ \ \ \
\ \ \ \ \ \ \ \ \ \ \ \ \ \ \ \ \ \ \ \ \ \ \ \ \ \ \ \ \ \ \ \ \ \ \ \ \ \ $%
\square $ \newline
\textbf{Acknowledgment. }The author M. A. A. Khan would like to acknowledge
the support provided by the Higher Education Commission of Pakistan for
funding this work through project No. NRPU 5332.

\end{document}